\documentclass{aptpub}
\usepackage{setspace}

\authornames{Daniel Clark} 
\shorttitle{Point process deconvolution}

\usepackage{setspace}

\usepackage{amsmath}
\usepackage{graphicx}

\usepackage{verbatim,amsmath}

\begin{document}
\title{Deconvolution of point processes}

\authorone{Daniel Edward Clark}
\address{School of Engineering and Physical Sciences,
Heriot-Watt University, Edinburgh, EH14 4AS UK}



\begin{abstract}
The superposition of two independent point processes can be described
by multiplication of their probability generating functionals (p.g.fl.s).
The inverse operation, which can be viewed as a deconvolution,
is defined by dividing the superposed process by one of its constituent 
p.g.fl.s.
The  deconvolved process
is computed using the higher-order chain rule for G\^{a}teaux differentials.
The higher-order quotient rule for 
G\^{a}teaux differentials
is first established and then applied 
to point processes.

\end{abstract}

\keywords{
point processes;
generating functionals;
G\^{a}teaux differentials
}

\ams{60G55}{49J50}

The superposition of two point processes
can be formed by multiplication of 
their probability generating functionals~\cite{Moyal62}.
While it is clear that the inverse operation would involve division of the probability generating functionals, there are no results regarding the Janossy densities after this operation.
This is probably due to the fact that the denominator
involves an infinite sum
and there was no easy way of computing the 
point process.
In this article, I derive the deconvolved point process using the recently derived higher-order 
chain rule for G\^{a}teaux differentials~\cite{hocr}.

The paper is structured as follows.
Section 1 describes the necessary background 
on functionals to derive the result,
including the probability generating functional (p.g.fl.),
the G\^{a}teaux differential,
Janossy densities and the higher-order chain rule.
Section 2 presents the main result by first 
 establishing the $n^{th}$-order 
G\^{a}teaux
differential of the reciprocal of a functional in Lemma 4.
Theorem 1 establishes the higher-order quotient
rule for G\^{a}teaux differentials
using Leibniz' rule and the
recently derived Faa di Bruno's formula for G\^{a}teaux differentials~\cite{hocr}.
This is used to determine the 
expression for the
deconvolved point process
via the probability generating functional~\cite{Moyal62}.


\section{Probability generating functionals}
This section discusses generating functionals
as a means of describing systems
with a variable number of objects.
 G\^{a}teaux differentials are discussed as a means 
of detemining their constituent functions.
The higher-order chain rule is introduced
to determine the differentials of composite 
functionals~\cite{hocr}.

\newpage
\noindent
{\bf Definition:  probability    generating functional}
\\
Following Moyal~\cite{Moyal62}, define the {\it probability generating functional}
(p.g.fl.) 
with
\begin{align}
    G_X(\psi)&
    =p_0+
\sum_{n=1}^\infty
    {1\over n!}
    \left(
    \int
    \prod_{i=1}^n
    dx_i\psi(x_i)
    \right)
    p_{n}(x_1,\ldots,x_n),
    \end{align}
where $p_{n}(x_1,\ldots,x_n)\ge0$ and
\begin{align}
p_0+
\sum_{n=1}^\infty
    {1\over n!}
    \left(
    \int
    \prod_{i=1}^n
    dx_i
    \right)
    p_{n}(x_1,\ldots,x_n)=1
\end{align}
The p.g.fl. is commonly used in point process theory
to characterise point patterns~\cite{RCoxIsham}.
For compactness of notation,
we write the p.g.fl. as 
\begin{align}
    G_X(\psi)&
				    =								\sum_{n=0}^\infty
								    {1\over n!}
												    \left(
																    \int
																				    \prod_{i=1}^n
																								    dx_i\psi(x_i)
																												    \right)
																    p_{n}(x_1,\ldots,x_n).
																				    \end{align}

\medskip
\noindent
{\bf Definition: G\^{a}teaux differential}
\\
A  {\it G\^{a}teaux differential}~\cite[p71-74]{Hille},
of a functional $\Psi(\psi)$
with increment $\xi$, can be defined with
\begin{align}
\delta\Psi(\psi;\xi)
=
\lim_{\epsilon\rightarrow 0}
{1\over\epsilon}
\left(\Psi(\psi+\epsilon\xi)-\Psi(\psi)\right).
\end{align}
The $n^{th}$-order differential (or variation), 
which is defined by recursively,
is denoted $\delta^n\Psi(\psi,\xi_1,\ldots,\xi_n)$.

\medskip
\noindent
{\bf Lemma 1}\\
The p.g.fl. $G(h)$ defines $P$
uniquely, since the Janossy densities can be recovered with
G\^{a}teaux differentials~\cite{Gateaux}, i.e.
\begin{align}
p_k(x_1,\ldots,x_k)=
\delta^k
G(0;
{\xi_1,\ldots,\xi_k}
),
\end{align}
where 
$\xi_i(x)=\delta_{x_i}(x)$.

\noindent
{\bf Proof}
\\
See Moyal~\cite{Moyal62}.

\medskip
\noindent
{\bf Lemma 2  }
\\
{
Let $\Pi$ be the set of all  partitions of
variables $\eta_1,\ldots,\eta_n$, and
$\pi\in\Pi$ denote a single partition
that has constituent blocks $\omega\in\pi$ of size $|\omega|$
consisting of constituent elements
$\zeta_{\pi,\omega,1},\ldots, \zeta_{\pi,\omega,|\omega|}\in \{\eta_1,\ldots,\eta_n\}$.
The $n^{th}$-order variation of composition $f\circ g$
with increments $\eta_1, \ldots, \eta_n$ is given by
}
\begin{align}
\label{ee}
\delta^n
\left(
\left(f\circ g\right)\left(y\right);
{\eta_1,\ldots,\eta_n}
\right)
=
\sum_{\pi\in \Pi}
\delta^{|\pi|}f
\left(
g(y);
\xi_{\pi,1}(y),\ldots,\xi_{\pi,|\pi|}(y)
\right),
\end{align}
{
where $\xi_{\pi,\omega}(y)$
is the variation of order $|\omega|$ with increments
$\zeta_{\pi,\omega,1},\ldots,\zeta_{\pi,\omega,{|\omega|}}$, i.e.}
\begin{align}
\xi_{\pi,\omega}(y)=
\delta^{|\omega|}g
\left(
y;
\zeta_{\pi,\omega,1},\ldots,\zeta_{\pi,\omega,{|\omega|}}
\right).
\end{align}

\noindent
{\bf Proof}
\\
See Clark~\cite{hocr}.

\medskip
\noindent
{\bf Lemma 3  }
\\
Leibniz' rule for functional derivatives 
the differentials of the product
$f(y)g(y)$ with increments
$\Xi=\{\eta_1,\ldots,\eta_n\}$
is
\begin{align}
\delta^n\left(f(y) {g(y)};\eta_1,\ldots,\eta_n\right)
=
\sum_{\Phi\subset \Xi}\delta^{|\Phi|}f\left(y; \phi_1,\ldots,\phi_{|\Phi|}\right)
\delta^{|\Xi|-|\Phi|}g\left(y; \xi_1,\ldots,\xi_{|\Xi|-|\Phi|}\right),
\end{align}
where the
summation is over all subsets $\Phi$ of the increments,
and $\Phi=\{\phi_1,\ldots,\phi_{|\Phi|}\}$, 
$\Xi-\Phi=\{\xi_1,\ldots,\xi_{|\Xi|-|\Phi|}\}$.

\noindent
{\bf Proof}
\\
See Mahler~\cite[p389]{Mahlerbook}.

\section{Deconvolution of point processes}
This section presents the main result.
Lemma 4 establishes the higher-order differential of the reciprocal
of a functional,
which is then used to establish a formula for the higher-order quotient
rule for G\^{a}teaux differentials in Theorem 1.
It is then shown in that the result is directly applicable to point processes.

\medskip
\noindent{\bf Lemma 4}\\
Consider a functional $g(y)$ in Banach space $\chi$.
Then the
$n^{th}$-order G\^{a}teaux differential~\cite{Gateaux} of $1/g(y)$
with set of increments $\Pi=\{\eta_1,\ldots,\eta_n\}$
is found with
\begin{align}
\delta^n(1/g(y);\eta_1,\ldots,\eta_n)=
\sum_{\pi\in\Pi}
{(-1)^{|\pi|}|\pi|!
\over
g(y)^{|\pi|+1}}
\prod_{\omega\in\pi}
\delta^{|\omega|}g
\left(
y;
\zeta_{\pi,\omega,1},\ldots,\zeta_{\pi,\omega,{|\omega|}}
\right),
\end{align}
where the summation
is over all partitions $\pi$ of 
increment set $\Pi$.
The product is over all cells
$\omega=\{\zeta_{\pi,\omega,1},\ldots,\zeta_{\pi,\omega,{|\omega|}}\}$ of partition $\pi$.

\medskip

\noindent{\bf Proof}\\
Suppose that $f(x)=1/x$, where $x$ is real-valued and non-zero.
Then, following the
definition of the G\^{a}teaux derivative (which in this case
reduces to the same rules as ordinary derivative),
\begin{align}
\delta^n
\left( 1/x; 
{\eta_1,\ldots,\eta_n}
\right)
=
{(-1)^nn!
\over
x^{n+1}
}.
\end{align}
Substituting this into (\ref{ee}) gives the required result.


\medskip

\noindent{\bf Theorem 1}\\
The higher-order quotient rule for G\^{a}teaux differentials
with increments
$\Xi=\{\eta_1,\ldots,\eta_n\}$
is determined with

\begin{align}
&\delta^n\left({f(y) \over g(y)};\eta_1,\ldots,\eta_n\right)
=
\\\notag
&
\sum_{\Pi\subset\Xi}
\left(
\sum_{\pi\in\Pi}
{(-1)^{|\pi|}|\pi|!
\over
g(y)^{|\pi|+1}}
\prod_{\omega\in\pi}
\delta^{|\omega|}g\left(y;\zeta_{\pi,\omega,1},\ldots,\zeta_{\pi,\omega,{|\omega|}}\right)
%
%
\right)
\delta^{|\Xi|-|\Pi|}f\left(y; \phi_1,\ldots,\phi_{|\Xi|-|\Pi|}\right),
\end{align}
where the summation on the left is
over all subsets $\Pi$ of $\Xi$.
The inner summation is
over all partitions $\pi$ of 
subset $\Pi$.
The term on the right,
$\delta^{|\omega|}g\left(y;\zeta_{\pi,\omega,1},\ldots,\zeta_{\pi,\omega,{|\omega|}}\right)$
is the variation of order $|\Xi|-|\Pi|$
with increments $\phi_1,\ldots,\phi_{|\Xi|-|\Pi|}$ that 
comprise $\Xi-\Pi$.


\medskip

\noindent{\bf Proof }\\
The result follows from application of Lemmas 3 and 4. 

\medskip
\noindent{\bf Corollary}\\
Consider a point process defined
by the superposition of point processes
with p.g.fl.s, i.e. 
\begin{align}
G[h]=G_1[h]G_2[h],
\end{align}
defined with Janossy
densities $p(x_1,\ldots,x_n)$, for $G[h]$,
$q(x_1,\ldots,x_n)$, for $G_1[h]$,
and $r(x_1,\ldots,x_n)$, for $G_2[h]$,
for $n>0$.
Let
$\Xi=\{x_1,\ldots,x_n\}$.
Then we can recover the Janossy densities
$r(x_1,\ldots,x_n)$
of $G_2[h]$
with
\begin{align}
&r(x_1,\ldots,x_n)
=
\\&
\sum_{\Pi\subset\Xi}
\left(
\sum_{\pi\in\Pi}
{(-1)^{|\pi|}|\pi|!
\over
q_0^{|\pi|+1}}
\prod_{\omega\in\pi}
q\left(v_{\pi,\omega,1},\ldots,v_{\pi,\omega,{|\omega|}}\right)
\right)
p\left(w_1,\ldots,w_{|\Xi|-|\Pi|}\right),
\end{align}
where the summation on the left is over all subsets $\Pi$ of
$\Xi$.
The inner summation is over all partitions $\pi$ of subset $\Pi$.
The term on the right $p\left(w_1,\ldots,w_{|\Xi|-|\Pi|}\right)$
is the Janossy density evaluated at 
$w_1,\ldots,w_{|\Xi|-|\Pi|}$
which comprise subset $\Xi-\Pi$,
and $q_0$ is the zeroth term in $G_1[h]$.
When $n=0$,
$r_0=p_0/q_0$.
%
\\

\medskip
\noindent{\bf Proof}
\\
Let $G_1[h]=f(y)$ and $G_2[h]=g(y)$
in (i).
Using  Lemma 1,
and setting $\eta_1=\delta_{x_1},\ldots,\eta_n=\delta_{x_n}$, and $h=0$ yields the result for $n>0$.
For $n=0$, we note that $G[0]=p_0$, and $G_1[0]=q_0$,
hence $r_0=p_0/q_0$.

\bibliography{bib_new1}
\end{document}